\documentclass[11pt,reqno, english]{amsart}

\usepackage{amssymb}
\usepackage{amscd}
\usepackage{amsfonts}
\usepackage{mathrsfs}
\usepackage{setspace}
\usepackage{version}
\usepackage[pdftex,colorlinks]{hyperref}

\usepackage{babel}
\usepackage[babel]{csquotes}

\numberwithin{equation}{section} \theoremstyle{plain}
\newtheorem{theorem}[subsection]{Theorem}
\newtheorem{thm-def}[subsection]{Theorem-Definition}
\newtheorem{proposition}[subsection]{Proposition}
\newtheorem{lemma}[subsection]{Lemma}

\newtheorem{definition}[subsection]{Definition}

\newcommand{\eps}{\varepsilon}

\newcommand\Z{\mathbb{Z}}
\newcommand\R{\mathbb{R}}

\newcommand\C{\mathbb{C}}
\newcommand\N{\mathbb{N}}

\newcommand\SL{\operatorname{SL}}

\newcommand\GL{\operatorname{GL}}

\newcommand\Q{\mathbb{Q}}

\newcommand{\sx}{{\mathsf x}}

\begin{document}

\title{On metric diophantine approximation in matrices and Lie groups}

\author{M. Aka, E. Breuillard, L. Rosenzweig and N. de Saxc\'e},


\address{Departement Mathematik\\
ETH Z\"urich\\
R\"amistrasse 101\\
8092 Zurich\\
SWITZERLAND}
\email{menashe-hai.akka@math.ethz.ch}

\address{Laboratoire de Math\'ematiques\\
B\^atiment 425, Universit\'e Paris Sud 11\\
91405 Orsay\\
FRANCE}
\email{emmanuel.breuillard@math.u-psud.fr}

\address{Department of mathematics\\
KTH\\
SE-100 44  Stockholm\\
SWEDEN}
\email{lior.rosenzweig@gmail.com}

\address{LAGA Institut Galil\'ee\\
Universit\'e Paris 13\\
93430 Villetaneuse\\
FRANCE}
\email{desaxce@math.univ-paris13.fr}

\begin{abstract}
We study the diophantine exponent of analytic submanifolds of $m \times n$ real matrices, answering questions of Beresnevich, Kleinbock and Margulis. We identify a family of algebraic obstructions to the extremality of such a submanifold, and give a formula for the exponent when the submanifold is algebraic and defined over $\Q$. We then apply these results to the determination of the diophantine exponent of rational nilpotent Lie groups.
\end{abstract}

\maketitle

\noindent {\bf Introduction}
\bigskip

In their breakthrough paper \cite{kleinbock-margulis}, Kleinbock and Margulis have solved a long standing conjecture of Sprindzuk regarding metric diophantine approximation on submanifolds of $\R^n$, stating roughly speaking that non-degenerate submanifolds are extremal in the sense that almost every point on them has similar diophantine properties as a random vector in $\R^n$ (i.e. it is not very well approximable, see below). Doing so they used new methods coming from dynamics and based on quantitative non-divergence estimates (going back to early work of Margulis \cite{margulis1971} and Dani \cite{dani1985}) for certain flows on the non-compact homogeneous space $\SL_n(\R)/\SL_n(\Z)$. They suggested at the end of their paper to extend their results to the case of submanifolds of matrices $M_{m,n}(\R)$, a natural set-up for such questions. This was studied further in \cite{kleinbock-margulis-wang}, \cite{beresnevich-kleinbock-margulis} and the problem appears in Gorodnik's list of open problems \cite{gorodnik-AIM}.

In this note we announce a set of results \cite{abrs2}, which give what we believe is a fairly complete picture of what happens in the matrix case as far as extremality is concerned. We identify a natural family of obstructions to extremality (Theorem \ref{extremal}) and show that they are in some sense the only obstructions to be considered (Theorem \ref{closure}). Our results also extend to the matrix case previous work of Kleinbock \cite{kleinbock-extremal-gafa,kleinbock-anextension} regarding degenerate submanifolds of $\R^n$. When the submanifold is algebraic and defined over $\Q$ we obtain a formula for the exponent (Theorem \ref{rational}).

In a second part of this note, we state new results regarding diophantine approximation on Lie groups, in the spirit of our earlier work \cite{abrs}. These results, which are applications of the theorems described in the first part of this note, concern the diophantine exponent of nilpotent Lie groups and were our initial motivation for studying diophantine approximation on submanifolds of matrices. The submanifolds to be considered here are images of certain word maps. Depending on the structure of the Lie algebra and its ideal of laws, these submanifolds can be degenerate. The relevant obstructions can nevertheless be identified and this leads to a formula for the diophantine exponent of an arbitrary rational nilpotent Lie group (Theorem \ref{formula}). A number of examples are also worked out explicitly.

\section{Diophantine approximation on submanifolds of $\R^n$.}   A vector $x \in \R^n$ is called \emph{extremal} (or \emph{not very well approximable}), if for every $\eps>0$ there is $c_\eps>0$ such that
$$|q \cdot x + p| >  \frac{c_\eps }{\|q\|^{n+ \eps}}  $$
for all $p \in \Z$ and all $q \in \Z^n \setminus\{0\}$. Here $q \cdot x$ denotes the standard scalar product in $\R^n$ and $\|q\|:=\sqrt{q \cdot q}$ the standard Euclidean norm.

As is well-known (Borel-Cantelli) Lebesgue almost every $x \in \R^n$ is extremal. An important question in metric diophantine approximation is that of understanding the diophantine properties of points $x$ that are allowed to vary inside a fixed submanifold $\mathcal{M}$ of $\R^n$. The submanifold $\mathcal{M}$ is called \emph{extremal} if Lebesgue almost every point on $\mathcal{M}$ is extremal. A key result here is
\bigskip

\begin{theorem}[Kleinbock-Margulis, \cite{kleinbock-margulis}]\label{km} Let $U$ be an open connected subset of $\R^k$ and $\mathcal{M}:=\{\mathbf{f}(x); x \in U\}$, where $\mathbf{f}:U \to \R^n$ is a real analytic map. Assume that $\mathcal{M}$ is not contained in a proper affine subspace of $\R^n$, then $\mathcal{M}$ is extremal.
\end{theorem}

\bigskip
This answered a conjecture of Sprindzuk. The proof made use of homogeneous dynamics via the so-called \emph{Dani correspondence} between diophantine exponents and the rate of escape to infinity of a diagonal flow in the space of lattices. We will also utilize these tools.

\section{Diophantine approximation on submanifolds of matrices.} It is natural to generalize this setting to that of submanifolds of matrices, namely submanifolds $\mathcal{M} \subset M_{m,n}(\R)$. The diophantine problem now becomes that of finding good integer approximations (by a vector $p \in \Z^m$) of the image $M \cdot q$ of an integer vector $q \in \Z^n$ under the linear endomorphism $M \in M_{m,n}(\R)$. The case $m=1$ corresponds to the above classical case (that of linear forms), while the dual case $n=1$ corresponds to simultaneous approximation. %The two are related by the transference principle (\cite[ch. 5]{cassels}).

% It is convenient to introduce the following definition:

% \begin{definition}[Diophantine exponent] We say that a matrix $M \in M_{m,n}(\R)$ has diophantine exponent $\beta(M) \geq 0$, if for $\beta(M)$ is the supremum of all numbers $\beta \geq 0$ for which there are infinitely many $(p,q) \in \Z^m \times \Z^n$, $q \neq 0$, such that
% $$ \|M \cdot q + p\| < \frac{1}{\|q\|^{\beta}}.$$
% \end{definition}

%The matrix $M$ is said to be \emph{diophantine} if $\beta$ is finite and \emph{extremal} if $\beta(M)=\frac{m}{n}$. Here the analogous pigeonhole or Dirichlet argument gives that $\beta(M) \geq \frac{n}{m}$ for every $M \in M_{m,n}(\R)$. Similarly, we say that the submanifold $\mathcal{M} \subset M_{m,n}(\R)$ is \emph{extremal} if Lebesgue almost every matrix in $\mathcal{M}$ is extremal.

%The problem of giving satisfying geometric conditions on $\mathcal{M}$ for it to be extremal, and the more general question of determining the exponent $\beta(M)$ of a typical matrix in $\mathcal{M}$ was asked by Kleinbock and Margulis at the end of their original paper \cite{kleinbock-margulis} and studied in the papers \cite{kleinbock-margulis-wang, beresnevich-kleinbock-margulis}. It is also mentioned as an open problem in Gorodnik's list \cite{gorodnik-AIM}.

%\bigskip

%In this note we prove Theorems \ref{extremal,general,closure,rational} below, which give an answer to these questions, and apply them to compute the diophantine exponent of nilpotent Lie groups, which was our original motivation (see Theorem \ref{liegroup} below).

%\bigskip

It turns out that it is more natural to study the slightly more general problem of approximating $0$ by the image $M \cdot q$ of an integer vector $q$. One can pass from the old problem to the new by embedding $\mathcal{M}$ inside $M_{m,m+n}(\R)$, via the embedding ($I_m$ denotes the $m \times m$ identity matrix)

\begin{eqnarray*}
M_{m,n}(\R) &\to& M_{m,m+n}(\R)\\
M &\mapsto& (I_m | M)
\end{eqnarray*}

From now on, we will consider an arbitrary connected analytic submanifold $\mathcal{M} \subset M_{m,m+n}(\R)$, given as $\mathcal{M}:=\{\mathbf{f}(x); x \in U\}$, where $\mathbf{f}:U \to M_{m,m+n}(\R)$ is a real analytic map from a connected open subset $U$ in some $\R^k$.

\bigskip

\begin{definition}[Diophantine exponent] We say that a matrix $M \in M_{m+n,n}(\R)$ has diophantine exponent $\beta(M) \geq 0$, if $\beta(M)$ is the supremum of all numbers $\beta \geq 0$ for which there are infinitely many $q \in \Z^{m+n}$ such that
 $$ \|M \cdot q\| < \frac{1}{\|q\|^{\beta}}.$$
 \end{definition}

\section{The pigeonhole argument and the obstructions to extremality.}\label{pigeon} By the pigeonhole principle (Dirichlet's theorem), the lower bound $\beta(M) \geq \frac{m}{n}$ holds for all $M$. Indeed one compares the number of integer points in a box of side length $T$ in $\Z^{m+n}$ with the volume occupied by the image of this box under $M$ in $\R^m$.
%Indeed up to multiplicative constants, there are roughly $T^{m+n}$ integers vectors $q \in \Z^{m+n}$ whose coordinates are less that $T$ in absolute value: under the linear map $M$, this gives rise to (a multi-set of) roughly $T^{m+n}$ points in $\R^m$, all having their coordinates bounded by a fixed multiple of $T$, hence occupying a volume of size roughly $T^m$ in $\R^m$. By the pigeonhole principle at least two of these points must lie less than $(T^{m+n}/T^m)^{-1/m}=T^{-n/m}$ apart, hence $\beta(M) \geq \frac{n}{m}$.
Furthermore, instead of considering the full box of side length $T$ in $\Z^{m+n}$, we could have restricted attention to the intersection of this box with a \emph{rational} subspace $W \leq \R^{m+n}$. The same argument would have then given the lower bound $$\beta(M) \geq \frac{\dim W}{\dim MW} -1.$$ Of course it may happen, given $M$, that for some exceptional subspace $W$, $\frac{\dim W}{\dim MW} -1 > \frac{n}{m} = \frac{n+m}{m}-1$. And this may well also happen for all $M \in \mathcal{M}$, provided $\mathcal{M}$ lies in the following algebraic subvariety $\mathcal{P}_{W,r}$ of $M_{m,m+n}(\R)$

\begin{equation}\label{pencil}\mathcal{P}_{W,r} : = \{ M \in M_{m,m+n}(\R); \dim MW \leq r\},
\end{equation}
where $W$ is a fixed rational subspace of $\R^{m+n}$ and $r$ a non-negative integer such that
\begin{equation}\label{ob}
\frac{\dim W}{r} -1 > \frac{n}{m}.
\end{equation}

By convention, we agree that $(\ref{ob})$ is satisfied if $r=0$. We will call the subvariety $\mathcal{P}_{W,r}$  of $M_{m,m+n}(\R)$ a \emph{pencil of endomorphisms} with parameters $W$ and $r$ (defined also for arbitrary, non rational, subspaces $W$). Note that when $m=1$, and $r=0$, this notion reduces to the notion of linear subspace (the orthogonal of $W$) of $\R^{n+1}$ (or affine subspace of $\R^n$). Hence asking that the submanifold $\mathcal{M}$ be not contained in any of those pencils $\mathcal{P}_{W,r}$ satisfying $(\ref{ob})$ is analogous in the matrix context to the condition of Theorem \ref{km} that $\mathcal{M}$ be not contained in an affine subspace.

\bigskip

\begin{theorem}[Extremal submanifolds]\label{extremal} Let $\mathcal{M} \subset M_{m,m+n}(\R)$ be a connected real analytic submanifold. Assume that $\mathcal{M}$ is not contained in any of the pencils $\mathcal{P}_{W,r}$, where $W,r$ range over all non-zero linear subspaces $W \leq \R^{m+n}$ and non-negative integers $r$ such that $(\ref{ob})$ holds. Then $\mathcal{M}$ is extremal, i.e. $\beta(M)=\frac{n}{m}$ for Lebesgue almost every $M \in \mathcal{M}$.
\end{theorem}

\bigskip

This result is close in spirit to that of \cite{beresnevich-kleinbock-margulis}, which gave a sufficient geometric  condition for strong extremality. The condition in Theorem \ref{extremal} is strictly weaker. It does not imply strong extremality, but only extremality, and with regards to extremality it is the optimal condition, as shown below in Theorems \ref{closure} and \ref{rational}.

\subsection{Non extremal submanifolds}
A general result of Kleinbock \cite{kleinbock-almostallversus} implies that the diophantine exponent of a random point of $\mathcal{M}$ is always well-defined. Namely there is $\beta= \beta(\mathcal{M})\in [0,+\infty]$ such that for Lebesgue almost every $x  \in  U$,

$$\beta(\mathbf{f}(x)) = \beta(\mathcal{M}).$$

Our first result is a general upper bound:

 \begin{theorem}[Upper bound on the exponent] \label{general} Let $\mathcal{M} \subset M_{m,m+n}(\R)$ be an analytic submanifold as defined above. Then
 $$\beta(\mathcal{M}) \leq \max \{ \frac{\dim W}{r} -1 ; \mathcal{P}_{W,r} \supset \mathcal{M} \}.$$
 \end{theorem}

Of course Theorem \ref{extremal} is an immediate consequence of this bound.

\bigskip

In \cite{kleinbock-extremal-gafa,kleinbock-anextension} Kleinbock showed that the diophantine exponent of an analytic submanifold of $\R^n$ depends only on its linear span. Our next result is a matrix analogue of this fact. Note that the diophantine exponent of a matrix $M$ depends only on its kernel $\ker M$. As $M$ varies in the  submanifold $\mathcal{M} \subset M_{m,m+n}(\R)$, consider the set of these kernels as a subset of the Grassmannian and take its linear span in the Pl\"ucker embedding. Denote by $\mathcal{H(M)}$ the matrices $M$ whose kernel lies in this linear span. The set $\mathcal{H(M)}$ is an algebraic subvariety containing $\mathcal{M}$ and contained in every pencil containing $\mathcal{M}$.

\bigskip

\begin{theorem}[Optimality of the exponent] \label{closure} We have:
$$\beta(\mathcal{M}) = \beta(\mathcal{H(M)}).$$ In particular $\beta(\mathcal{M})= \beta(\textnormal{Zar}(\mathcal{M}))$, where $\textnormal{Zar}(\mathcal{M})$ denotes the Zariski closure of $\mathcal{M}$, and $\beta(\mathcal{M}) = \beta(\Omega)$ for any open subset $\Omega \subset \mathcal{M}$.
\end{theorem}
\bigskip
In particular $\mathcal{M}$ is extremal if and only if $\mathcal{H(M)}$ is extremal.

\section{Lower bounds on the exponent and rationality}
Theorem \ref{general} gives a general upper bound on the exponent. The pigeonhole argument described at the beginning of \S \ref{pigeon} yields a lower bound on $\beta(\mathcal{M})$ in terms of the exponents of the \emph{rational obstructions} in which $\mathcal{M}$ is contained, i.e. the pencils $\mathcal{P}_{W,r}$ with $W$ a rational subspace of $\R^{m+n}$. Hence, for a general analytic submanifold $\mathcal{M} \subset M_{m,m+n}(\R)$, we only have the following general upper and lower bound:

\begin{equation}\label{uplow}
\max_{\mathcal{P}_{W,r} \supset \mathcal{M} , W \textnormal{ rational } } \frac{\dim W -r}{r}  \leq \beta(\mathcal{M}) \leq \max_{\mathcal{P}_{W,r} \supset \mathcal{M}} \frac{\dim W -r}{r}.
\end{equation}

For a submanifold $\mathcal{M}$ in general position the upper and lower bound are typically distinct. However  we will prove:

\bigskip

\begin{theorem}[Subvarieties defined over $\Q$] \label{rational} Assume that the Zariski-closure of the connected real analytic submanifold $\mathcal{M} \subset M_{m,m+n}(\R)$ is defined over $\Q$. Then the upper and lower bounds in $(\ref{uplow})$ coincide, and hence are equal to $\beta(\mathcal{M})$. In particular $\beta(\mathcal{M}) \in \Q$.
\end{theorem}

\bigskip
%Recall that the Zariski-closure of $\mathcal{M}$ is defined over $\Q$ if and only if the set of polynomial maps vanishing on $\mathcal{M}$ is invariant under the action of the Galois group $Gal(\R|\Q)$, or equivalently if it is defined by the vanishing of finitely many polynomial maps on $M_{m,m+n}(\R)$ with rational coefficients.

%One verifies that  $\beta(\mathcal{P}_{W,r}) \leq \max\{ \frac{\dim W -r}{r} , \frac{n}{m}\}$, for all $W$, and that equality holds if $W$ is rational. Hence %in the case when $\mathcal{H(M)}$ is defined over $\Q$, one has:
%$$\beta(\mathcal{M})= \beta(\mathcal{H(M)})= \max_{\mathcal{P}_{W,r} \supset \mathcal{M}} \beta(\mathcal{P}_{W,r}).$$

The proof of Theorem \ref{rational} is based on the following combinatorial lemma, which is used here with $G=\textnormal{Gal}(\C|\Q)$ and will be used once again later on in the applications to nilpotent groups with $G=\GL_k$.
%\bigskip

Let $V$ be a finite dimensional vector space over a field and $\phi: \textnormal{Grass}(V) \to \N \cup\{0\}$ a function on the Grassmannian, which is non-decreasing (for set inclusion) and \emph{submodular} in the sense that for every two subspaces $W_1$ and $W_2$ we have
$$\phi(W_1 + W_2) + \phi(W_1 \cap W_2) \leq \phi(W_1) + \phi(W_2).$$

\bigskip
\begin{lemma}[Submodularity lemma]\label{submodular} Let $G$ be a group acting by linear automorphisms on $V$. If $\phi$ is invariant under $G$, then the following minimum is attained on a $G$-invariant subspace
$$\min_{W \in \textnormal{Grass}(V) \setminus\{0\}} \frac{\phi(W)}{\dim W}.$$
\end{lemma}

\section{Diophantine approximation on Lie groups}
Inspired by work of Gamburd-Jakobson-Sarnak \cite{gamburd-jakobson-sarnak} and Bourgain-Gamburd \cite{bourgain-gamburd} on the spectral gap problem for finitely generated subgroups of compact Lie groups, we defined in a previous article \cite{abrs} the notion of diophantine subgroup of an arbitrary Lie group $G$. The definition is as follows. Any finite symmetric subset $S:=\{1,s_1^{\pm 1}, \ldots,s_k^{\pm 1}\}$ in $G$ generates a subgroup $\Gamma \leq G$. If for all $n \in \N$
$$ \inf\{ d(1,\gamma) ; \gamma \in S^n \setminus\{1\}\} > \frac{1}{|S^n|^\beta},$$
then we say that $(\Gamma,S)$ is \emph{$\beta$-diophantine}. And we say that $\Gamma$ is \emph{diophantine} if it is $\beta$-diophantine for some finite $\beta$.
Here $d(\cdot, \cdot)$ denotes a fixed Riemannian metric on $G$ and $|S^n|$ is the cardinality of the $n$-th product set $S^n:=S \cdot \ldots \cdot S$. It is easily seen that being diophantine does not depend on the choice of $S$ or  $d(\cdot, \cdot)$. And if $G$ is nilpotent this is also true of being $\beta$-diophantine.

\bigskip
The connected Lie group $G$ is said to be \emph{diophantine on $k$ letters} if for almost every choice  of $k$ group elements $s_1,\ldots,s_k$ chosen independently with respect to the Haar measure, the subgroup they generate is diophantine. Finally one says that $G$ is \emph{diophantine} if it is diophantine on $k$ letters for every integer $k$.

\bigskip

While it is conjectured that all semisimple Lie groups are diophantine, there are examples of non-diophantine Lie groups. Indeed a construction was given in \cite{abrs} for each integer $k \in \N$ of a connected Lie group which is diophantine on $k$ letters, but not on $k+1$ letters. Our examples are certain nilpotent Lie groups without a rational structure. We showed in that paper that the first examples arise in nilpotency class $6$ and higher. In fact every nilpotent Lie group $G$ with nilpotency class at most $5$, or derived length at most $2$ (i.e. metabelian), is diophantine.

\section{Diophantine exponent of nilpotent Lie groups}

If $G$ is nilpotent, $|S^n|$ grows like $n^{\alpha_S}$, where $\alpha_S$ is an integer given by the Bass-Guivarc'h formula. If the $k$ elements $s_i$'s forming $S$ are chosen at random with respect to Haar measure, then $\alpha_S$ is almost surely a fixed integer, which is a polynomial in $k$ (see \cite{abrs}).
\bigskip

\begin{proposition}[Zero-one law]\label{zerone} Let $G$ be a simply connected nilpotent Lie group, and pick an integer $k \geq \dim G/[G,G]$. There is a number $\beta_k \in [0,+\infty]$, such that if $\beta> \beta_k$ (resp. $\beta< \beta_k$), then with respect to Haar measure almost every (resp. almost no) $k$-tuple in $G$ generates a $\beta$-Diophantine subgroup.
\end{proposition}
\bigskip
The proof of this is based on the ergodicity of the group of rational automorphisms of the free Lie algebra on $k$ letters acting on $(\textnormal{Lie}(G))^k$.  When the nilpotent Lie group $G$ is rational (i.e. admits a $\Q$-structure) the exponent $\beta_k$ can be computed explicitly using Theorem \ref{rational}. We have:

\bigskip

\begin{theorem}[A formula for the exponent]\label{formula} Assume that $G$ is a rational simply connected nilpotent Lie group. There is a rational function $F \in \Q(X)$ with coefficients in $\Q$ such that for all  large enough $k$,
$$\beta_{k} = F(k).$$
In particular $\beta_k \in \Q$. When $k \to \infty$, $\beta_k$ converges to a limit $\beta_{\infty}$ with $0<  \beta_{\infty} \leq 1$.
\end{theorem}

\bigskip
For example, if $G$ is the $(2m+1)$-dimensional Heisenberg group and $k \geq 2m$, then $\beta_k=1-\frac{1}{k}-\frac{2}{k^2}$. More generally if $G$ is any $2$-step nilpotent group not necessarily rational, then $\beta_k = (1-\frac{1}{k})\frac{1}{\dim[G,G]} - \frac{2}{k^2}$ for $k \geq \dim G/[G,G]$.

\bigskip

We also obtain closed formulas for $\beta_k$ in the case when $G$ is the group of $n \times n$ unipotent upper-triangular matrices, e.g. if $n=4$, and $k \geq 3$, then $\beta_k= \frac{k^3-k-3}{k^3+k^2-k}$. And in the case when $G$ is an $s$-step free nilpotent group on $m$ generators, e.g. if $m=2$ and $s=3$, then $\beta_k=\frac{k^3-k-6}{2(k^3+k^2-k)}$. These formulas involve the dimensions of the maximal (for the natural partial order on Young diagrams) irreducible $\GL_k$-submodule of the free Lie algebra on $k$ generators modulo the ideals of laws of $G$.

\bigskip

The reduction to Theorem \ref{rational} proceeds as follows. Since $k$ is large, one can restrict attention to the last term $G^{(s)}$ in the central descending series. Given a $\Z$-basis $e_1,\ldots,e_{m+n}$ of the $s$-homogeneous part of the relatively free Lie algebra of $G$ on $k$ generators $\mathcal{F}_{k,G}$ (see \cite{abrs}), the submanifold $\mathcal{M}_{k,G}$ of matrices to be considered is the image of $(\textnormal{Lie}(G))^k$ under the (polynomial) map sending $\sx \in (\textnormal{Lie}(G))^k$ to the $(n+m) \times m$ matrix whose columns are the $e_i(\sx)$. Here $m=\dim G^{(s)}$. Computing the exponent amounts to first identify the pencils $\mathcal{P}_{W,r}$ in which $\mathcal{M}_{k,G}$ sits and then compute the maximum of the ratios $\frac{\dim W}{r}$. Using the submodularity lemma (Lemma \ref{submodular}) applied for the $\GL_k$ action of linear substitutions we may restrict attention to those pencils corresponding to subspaces $W$ of $\mathcal{F}_{k,G}$ that are fully invariant ideals. Determining those ideals is usually possible, depending on $G$, thanks to the known representation theory of the free Lie algebra viewed as a $\GL_k$-module.

% etc, etc

% \appendix

% \section{}
% \label{}

% The Acknowledgements are an un-numbered section
%\section*{Acknowledgements}
% Acknowledgements text here

\bibliographystyle{alpha}
\bibliography{bibli}

%\begin{thebibliography}{00}
% please try to use the bibitem system -
% the references should be in alphabetical order of authors' names.
% Articles with a single author first, author will 1 co-author next,
% then author with several co-authors;


%\bibitem{label}

%\end{thebibliography}

\end{document}